\documentclass[12pt,]{article}
\usepackage[normal]{subfigure}
\usepackage{graphicx}
\usepackage{amssymb}
\usepackage{amscd}
\usepackage{amsmath}
\usepackage{amsfonts}
\usepackage{amsbsy}
\usepackage{fancyhdr}

\usepackage{epsfig}

\numberwithin{equation}{section}

\numberwithin{theorem}{section}
 
 \numberwithin{lemma}{section}
\date{}

\begin{document}
\date{}
\title{On the number of limit cycles of  polynomial Li\'enard systems\thanks{The project supported  by  the  National
Natural Science Foundation of China (10971139) and the Slovenian
Research Agency.} }

\author{ Maoan Han$^{a,b}$\thanks {Corresponding author. {\em E-mail address}: mahan@shnu.edu.cn; tel. +8621-64328672; fax +8621-64328672
},
  Valery G.  Romanovski$^{b,c}$\thanks {{\em E-mail
address}: valery.romanovsky@uni-mb.si.}\\
\small \emph{$^a$ Department of
Mathematics, Shanghai Normal University,}\\
\small \emph{ Shanghai
200234, P. R. China } \\
\small \emph{ $^b$ Faculty of Natural Science and Mathematics,}\\
\small \emph{University of Maribor,  SI-2000 Maribor, Slovenia}\\
\small \emph{$^c$ CAMTP - Center for Applied Mathematics and Theoretical
Physics, } \\
\small \emph{ University of Maribor, SI-2000 Maribor, Slovenia}}

\date{}
\maketitle

 {\bf Abstract:}\small  {  Li\'enard systems are very important  mathematical  models describing oscillatory processes arising in applied sciences.
In this paper,
we study polynomial Li\'enard systems of arbitrary degree on the plane, and develop
 a new method to obtain a lower bound of the maximal
number of limit cycles. Using  the method and basing  on some known results for lower degree we obtain  new  estimations of the number of  limit  cycles  in the systems which greatly improve
existing results. }


\hspace{3mm}

\section{Introduction  and main results}

Consider a polynomial Li\'enard system of the form
\begin{equation}\label{01}
\dot{x}=y,\quad \dot{y}=-g(x)-\varepsilon f(x)y,
\end{equation}
where $\varepsilon$ is a small parameter, $f(x)$, $g(x)$ are
polynomials in  $x$ of  degree $n$ and $m$,  respectively. The above
system is called a Li\'{e}nard system. It describes the dynamics of systems of one degree of
freedom under existence of a linear restoring force and a nonlinear dumping.   It was shown
by  Li\'{e}nard \cite{lien} that under some conditions on the functions $f(x)$ and $g(x)$  in the system arise auto-oscillations.
In the first half of the last century models based on the Li\'{e}nard system
were important for the development of  radio and vacuum tube technology. Nowadays the system
is widely used to describe oscillatory processes arising in various studies of mathematical  models of physical, biological, chemical, epidemiological, physiological,   economical and  many other phenomena (see e.g. \cite{Gla,Lib_sur} and references therein).

Our study is devoted to finding  Li\'{e}nard  systems which
admit   not a single, but few auto-oscillatory regimes   (limit cycles).
 Let $H(n,m)$ denote the
maximal number of limit cycles of system (\ref{01})
 on the plane for $\varepsilon$ sufficiently small. The
lower bound of $H(n,m)$  for the Li\'{e}nard system  has been widely studied. For general $m$
and $n$ using the averaging theory of order 3, the authors of  \cite{Lli}
gave the estimation
$$H(n,m)\geq[\frac{n+m-1}2],$$
 which generalized the earlier bound  $H(n,1)\geq
\left[\frac{n}{2}\right]$ of  Blows and Lloyd \cite{blows}.
 Recently, Han, Tian and Yu \cite{hty}  obtained the following improvement
$$H(n,m)\geq \max\{[\frac{m-2}{3}]+[\frac{2n+1}{3}],
[\frac{n-2}{3}]+[\frac{2m+1}{3}]\}$$ for $m,n\geq 2$. For $m=2$,
Han \cite{han1999} proved $
{H}(n,2)\geq\left[\frac{2n+1}{3}\right],\ n\geq 2.
 $
For $m=3$, Dumortier and Li \cite{dl} obtained $H(2,3)\geq 5$,
Christopher and Lynch \cite{ls} proved
$${H}(n,3)\geq 2\left[\frac{3n+6}{8}\right],\ \  2\leq n\leq
50,$$ and Yang, Han and Romanovski \cite{yhr} further obtained
$$H(n,3)\geq \left[\frac{3n+14}{4}\right], \ 2\leq n\leq 8,
$$
 which gave a larger lower bound of $H(n,3)$ than the above ones for $3\leq n\leq 8$.

From  \cite{hzy,yha} we know that
$$
H(n,3)\geq n+2-\left[\frac{n+1}{4}\right], \ \ 9\leq n\leq 22.
$$

For $m=4$,   Han, Yan, Yang and Lhotka \cite{hyy1} studied the limit
cycle bifurcation  of system \eqref{01}  and  obtained
$$H(n,4)\geq n+3,\ n=2,3,5,6,7,8, \ H(4,4)\geq 6.$$
Christopher and Lynch \cite{ls} gave $H(9,4)\geq 9.$ And then, Yu
and Han  \cite{yh,yhan} obtained $$H(n,4)\geq n,\ n=10, 11, 12, 13,
14.$$
$$H(n,4)\geq
n+4-\left[\frac{n+1}{5}\right],\ 3\leq n\leq 18.$$

More results for some concrete $m$ and $n$ can be found in
\cite{yh}.

Motivated by \cite{hl}, in this paper we give a new method to find a
lower bound of $H(n,m)$ for many integers $m$ and $n$. The main
results are the following.

$(1)$ \ $H(n,4)\geq H(n,3)\geq 2[\frac{n-1}{4}]+[\frac{n-1}{2}]$, $n\geq 3.$

$(2)$ \ $H(n,6)\geq H(n,5)\geq 2[\frac{n-1}{3}]+[\frac{n-1}{2}]$, $n\geq 5.$

$(3)$ \ $H(n,7)\geq\frac{3}{2}n-9$   for $ n\geq 7$. In general, for any integer $m\geq 7$ there exists $\gamma_m>0$ satisfying
$$\lim_{m\rightarrow\infty} \sup \frac{\gamma_m}{(m+1)\ln(m+1)}\leq \frac{1}{4\ln2}$$
such that $$H(n,m)\geq (\frac{\ln(m+2)}{2\ln2}-\frac{1}{3})n-\gamma_m, \ \ n\geq m.$$ In particular,
$$\displaystyle\lim_{m\rightarrow\infty}\lim_{n\rightarrow\infty} \inf \frac{H(n,m)}{n\ln(m+2)}\geq
\frac{1}{2\ln2}.$$

$(4)$ \ For any integer $r\geq 0$, $$\lim_{m\rightarrow\infty}\inf\frac{H(m\pm r,m)}{m\ln m}\geq \frac{1}{2\ln2}.$$

$(5)$\ For $m=2^{p+1}-1$, $p\geq 1$, we have $$H(m-1,m)\geq \frac{(m+1)\ln(m+1)}{2\ln2}+1,$$and $$\label{44}H(m,m)\geq \frac{(m+1)\ln(m+1)}{2\ln2}+1.$$

$(6)$ \ For all $m\geq 3$ we have $$H(m-1,m)\geq \frac{(m+2)\ln(m+2)}{3\ln2}-\frac{m+2}{3}(1+\frac{\ln3}{\ln2})+1,$$ and

$$H(m,m)\geq \frac{(m+2)\ln(m+2)}{3\ln2}-\frac{m+2}{3}(1+\frac{\ln3}{\ln2})+1.$$

$(7)$\ Let $r$ be a positive integer. For any $k\geq 2$ there exist constants $B_{k,r}$ and $\bar B_{k,r}$ satisfying
$$ B_{k,r}\geq-\Big[\frac{1}{k}\Big(1+[\frac{r}{2}]\Big)+\frac{\ln k}{2\ln2}\Big],$$
$$\bar B_{k,r}\geq-\Big[\frac{1}{k+1}\Big(2+[\frac{r}{2}]\Big)+\frac{\ln(k+1)}{2\ln2}\Big],$$
such that
$$H(m-r,m)\geq\frac{(m+1)\ln(m+1)}{2\ln2}+B_{k,r}(m+1)+1+[\frac{r}{2}]$$
for $m=2^pk-1$, $ p\geq 1$, and
$$H(m-r,m)\geq\frac{(m+2)\ln(m+2)}{2\ln2}+\bar B_{k,r}(m+2)+2+[\frac{r}{2}]$$
for $m=2^p(k+1)-2$, $ p\geq 1$.

The conclusions listed above are contained in  Theorems 3.1, 3.2, 4.2, 4.3, 5.1, 5.2 and formula \eqref{50} in sections 3, 4 and 5 below.

\section{ Preliminary results }

The polynomial Li\'enard system \eqref{01} can be transformed  to the form

\begin{equation}\label{11}
\begin{array}{l}
\dot{x}=y-\varepsilon F(x),\ \dot{y}=-g(x),
\end{array}
\end{equation}
where $F(x)=\int_0^xf(x)dx.$

We introduce the following definition.

{\bf Definition 2.1.} We say that the system \eqref{11} has property
$Z(n,m,k)$ if the following are satisfied:

(1) $\deg F\leq n+1$, $\deg g\leq m$, and $\lim_{x\rightarrow
+\infty}g(x)=+\infty$;

(2) there are a constant $\varepsilon_0>0$ and a compact set
$D\subset \mathbb{R}^2$ such that for all $0<\varepsilon\leq \varepsilon_0$
the system \eqref{11} has at least $k$ limit cycles in $D$, each
having an odd multiplicity.

Obviously, if we can find a polynomial system \eqref{11} which has
the property $Z(n,m,k)$, then $H(n,m)\geq k$.

As we know, the first order Melnikov function of \eqref{11} has the
form
$$M(h)=\oint_{L_h}F(x)dy,$$
where $L_h$ is a smooth closed curve defined by the equation
$H(x,y)=h$
 on the plane, with
$$H(x,y)=\frac{1}{2}y^2+G(x), \ \ G(x)=\int_0^xg(x)dx.$$

As it is known \cite{LiJibin},
 \eqref{11} has the property $Z(n,m,k)$ if the function $M(h)$
  has at least $k$ zeros, each having an odd multiplicity.

From the work of \cite{blows,han1999,yhr,yhan}, we know that the
following facts hold.

{\bf Lemma 2.1.} \ { \it There are  polynomial Li\'enard systems of
the form \eqref{11} which have the following properties

$(1)$ $Z(n,1,[\frac{n}{2}])$ and $Z(n,2,[\frac{2n+1}{3}])$ for
$n\geq 1$;

$(2)$ $Z(n,3,[\frac{3n+14}{4}])$ for $2\leq n\leq 8$ and
$Z(n,4,n+4-[\frac{n+1}{5}])$ for $3\leq n\leq 18$.}

 Next we give a method to construct polynomial  Li\'enard systems
 having new properties starting from a given polynomial Li\'enard system
 having a certain property.  We begin with  the following lemma.

 {\bf Lemma 2.2.} \ {\it Let \eqref{11} have
the property $Z(n,m,k)$. Then there exists $x^*>0$ such that for all
$x_0<-x^*$ the system
\begin{equation}\label{21}
\begin{array}{l}
\dot{x}=y-\varepsilon F(x^2+x_0),\ \dot{y}=-2xg(x^2+x_0)
\end{array}
\end{equation}
has  the property  $Z(2n+1,2m+1,2k)$.}

{\bf Proof.} By Definition 2.1, there are a constant
$\varepsilon_0>0$ and a compact set $D\subset \mathbb{R}^2$ such that for all
$0<\varepsilon\leq \varepsilon_0$ the system \eqref{11} has at least
$k$ limit cycles in $D$, each having an odd multiplicity. Since $D$
is compact, for some $x^*>0$ we must have $D\subset \{(x,y)| \
|x|<x^*\}$. Let $x_0<-x^*$. Then for $|x|<x^*$ we have $x-x_0>0$.
Note that the change $u=x-x_0$ carries \eqref{11} into
\begin{equation}\label{22}
\begin{array}{l}
\dot{u}=y-\varepsilon F(u+x_0),\ \dot{y}=-g(u+x_0).
\end{array}
\end{equation}
Thus, the above system has $k$ limit cycles on a compact set which
is contained in $\{(u,y)| \ u>0\}$.

Further, we introduce $u=v^2$ to \eqref{22}  to obtain
$$ 2v\dot v=y-\varepsilon F(v^2+x_0), \ \ 2v\dot y=-2vg(v^2+x_0)$$
which is equivalent to
\begin{equation}\label{23}
\begin{array}{l}
\dot{v}=y-\varepsilon F(v^2+x_0),\ \dot{y}=-2vg(v^2+x_0)
\end{array}
\end{equation}
on $v<0$ or $v>0$. Therefore, \eqref{23} has $k$ limit cycles both
on a compact set $D_1$ in $v<0$ and on a compact set $D_2$ in $v>0$.
That is to say, it has $2k$ limit cycles on $D_1\bigcup D_2$. The
conclusion follows. This ends the proof.

{\bf Lemma 2.3.} \ {\it Suppose $G_1(x)$ is a polynomial of degree
$2l$ satisfying
$$G_1(x)=G_0x^{2l}(1+O(x^{-1})), \ l\geq 1, \ G_0>0$$
as $|x|\rightarrow \infty$. Let $H_1(x,y)=\frac{1}{2}y^2+G_1(x)$ and
$$I_j(h)=\oint_{H_1=h}x^{2j+1}dy, \ \ j\geq 0$$
along the orbits of the system
$$\dot x=y, \ \ \dot y=-G_1'(x).$$ Then there exist constants
$\alpha_j>0$, $j\geq 0$ such that
$$I_j(h)=-\alpha_jh^{\frac{2j+l+1}{2l}}(1+o(1))$$
as $h\rightarrow \infty$.}

{\bf Proof.} Noting $ydy=-G_1'(x)dx$ along $H_1(x,y)=h$, we have
$$
I_j(h)=-\displaystyle\oint_{H_1=h}\frac{x^{2j+1}G_1'(x)}{y}dx
=-2\displaystyle\int_{a(h)}^{b(h)}\frac{x^{2j+1}G_1'(x)}{\sqrt{2(h-G_1(x))}}dx,
$$
where $a(h)$ and $b(h)$ are the solutions of the equation $G_1(x)=h$
satisfying
$$\lim_{h\rightarrow\infty}a(h)=-\infty, \
\lim_{h\rightarrow\infty}b(h)=\infty.$$

Let $x_0>0$ and $x_0'<0 $ be such that
$$G_1(x_0)=G_1(x_0')\equiv u_0>0, \ G_1(x)>0, \ xG_1'(x)>0 \ {\rm
for} \ x_0'\leq x\leq x_0.$$ Then
$$I_j(h)=\tilde I_1(h)+\tilde I_2(h)+\tilde I_3(h),$$
where
\begin{equation}\label{24}
\tilde
I_1(h)=-\displaystyle2\int_{a(h)}^{x_0'}\frac{x^{2j+1}G_1'(x)}{\sqrt{2(h-G_1(x))}}dx=
2\int_{u_0}^{h}\frac{(x_1(h))^{2j+1}}{\sqrt{2(h-u)}}du,
\end{equation}
$$
\tilde
I_2(h)=-\displaystyle2\int_{x_0)}^{b(h)}\frac{x^{2j+1}G_1'(x)}{\sqrt{2(h-G_1(x))}}dx=
-2\int_{u_0}^{h}\frac{(x_2(h))^{2j+1}}{\sqrt{2(h-u)}}du,
$$
$$\tilde
I_3(h)=-\displaystyle2\int_{x_0'}^{x_0}\frac{x^{2j+1}G_1'(x)}{\sqrt{2(h-G_1(x))}}dx,$$
and  $x_1(u)<0<x_2(u)$ satisfy $G_1(x_i(u))=u$ for $u_0\leq u \leq
h$, $ i=1,2$. Obviously, $\lim_{h\rightarrow\infty}\tilde I_3(h)=0$.
Thus, to finish the proof we need only to prove that there exist
constants $\beta_1>0$ and $\beta_2>0$ such that
\begin{equation}\label{25}\tilde I_i(h)=-\beta_ih^{\frac{2j+l+1}{2l}}(1+o(1)),
\ i=1,2\end{equation} as $h\rightarrow \infty$. We only consider the
case $i=1$. The case $i=2$ is just similar.
 By our assumption on $G_1$, the equation $G_1(x)=u$ can be
 rewritten as
 $$|x|(1+O(|x|^{-1}))=(u/G_0)^{\frac{1}{2l}}$$
 for $u>0$ large. It then follows that
 $$x_1(u)=-(u/G_0)^{\frac{1}{2l}}(1+O(u^{\frac{-1}{2l}})).$$
 Hence
 $$\begin{array}{ll}
(x_1(h))^{2j+1}&=-\displaystyle(u/G_0)^{\frac{2j+1}{2l}}(1+O(u^{\frac{-1}{2l}}))\\
&=-\displaystyle(u/G_0)^{\frac{2j+1}{2l}}+u^{\frac{j}{l}}\varphi_1(u),\end{array}$$
where $\varphi_1(u)$ is smooth and bounded on the interval
$[u_0,+\infty)$.

Now introducing the  change $u=h\sin^2\theta$ we obtain from
\eqref{24}

$$\begin{array}{ll}
\tilde I_1(h))&=-\displaystyle2\int_{\theta_0(h)}^{\frac{\pi}{2}}\sqrt{2h}\sin \theta (x_1(h\sin^2\theta))^{2j+1}d\theta
\\
&=-\displaystyle h^{\frac{2j+1}{2l}+\frac{1}{2}}\tilde
\varphi_0(h)+h^{\frac{j}{l}+\frac{1}{2}}\tilde
\varphi_1(h),\end{array}$$ where
$$\tilde
\varphi_0(h)=\displaystyle2\sqrt{2}G_0^{-\frac{2j+1}{2l}}\int_{\theta_0(h)}^{\frac{\pi}{2}}[\sin
\theta]^{\frac{2j+1}{l}+1}d\theta,$$ $\tilde \varphi_1(h)$ is
bounded on $[u_0,+\infty)$, and $\theta_0(h)$ satisfies $\sin
\theta_0=\sqrt{u_0/h}$. It is evident that
$$\lim_{h\rightarrow\infty}\tilde
\varphi_0(h)=2\sqrt{2}G_0^{-\frac{2j+1}{2l}}\int_{0}^{\frac{\pi}{2}}[\sin
\theta]^{\frac{2j+1}{l}+1}d\theta\equiv\beta_1.$$

Then  \eqref{25} follows for $i=1$. This completes the proof.

Using the above lemma, we have further the following fundamental lemma.

 {\bf Lemma 2.4.} \ {\it Let
$G_2(x)$ and $F_2(x)$ be even polynomials in $x$ with
$G_2(\infty)=+\infty.$ Then for any integer $q\geq 1$ and a compact
set $U_0\subset \mathbb{R}^2$ there exist $\varepsilon_0>0$, $b_j\neq 0$,
$j=0,\cdots,q$ and a compact set $U\subset \mathbb{R}^2$ with $U\bigcap
U_0=\emptyset$ such that the following system

 \begin{equation}\label{26}
 \begin{array}{l}
\dot x=y-\displaystyle [\lambda F_2(x)+\mu
\sum_{j=0}^qb_jx^{2j+1}],\\
\dot y=-G_2'(x)\end{array}\end{equation} has $q$ limit cycles with
odd multiplicity in $U$ for all $|\lambda|\leq \varepsilon_0$ and
$0<|\mu|\leq \varepsilon_0$.}

{\bf Proof.} Since $G_2$ is even with $G_2(\infty)=+\infty$, there
exists $x_0>0$ such that for all $a\geq x_0$ the orbit $\gamma(a)$
of the Hamiltonian system $\dot x=y$, $\dot y=-G_2'(x)$ starting
from $A(a,0)$ is periodic, which ensures $G_2'(x)>0$ for $a>x_0$.
Also, for any given compact set $U_0\subset \mathbb{R}^2$ there is
$h^*>G_2(x_0)\equiv h_0$ such that $U_0\subset \{(x,y)| \
H_2(x,y)<h^*\}$.

 Hence, for any given $x_1>x_0$ satisfying $h^*<G_2(x_1)\equiv h_1$
there exists $\varepsilon_1=\varepsilon_1(x_1)>0$ such that for all
$a\in [x_0,x_1]$ and $|\lambda|\leq \varepsilon_1$ the orbit
$\gamma_\lambda(a)$ of the symmetric system $$\dot x=y-\lambda
F_2(x),\ \ \dot y=-G_2'(x)$$ starting from the same point $A(a,0)$
is also periodic. Then further, for any given $N>0$ there exists
$\varepsilon_2=\varepsilon_2(x_1,N)>0$ such that for all $a\in
[x_0,x_1]$, $|\lambda|\leq \varepsilon_1$, $|\mu|\leq \varepsilon_2$
and $|b_j|\leq N$, $j=0,\dots,q$ the system \eqref{26} has a
positive orbit $\gamma^+_{\lambda,\mu}(a)$  starting from the point
$A(a,0)$ which intersects the positive $x$-axis again at some point
$B(b(a,\lambda,\mu),0)$ for the first time. To find $\varepsilon_0$
and $U$,  suitable $x_1$ and $N$ will be chosen later.
 Let $$H_2(x,y)=\frac{1}{2}y^2+G_2(x).$$ Then along the orbit $\widehat{AB}$ of
 \eqref{26} we have
 $$\begin{array}{ll}
 H_2(B)-H_2(A)&=\displaystyle\int_{\widehat{AB}}(G_2'(x)dx+ydy)\\
 &=\displaystyle\int_{\widehat{AB}}(\lambda
 F_2(x)+\mu\sum_{j=0}^qb_jx^{2j+1})G_2'(x)dt.\end{array}$$

 Note that $b(a,\lambda,0)=a$ for all $a\in
[x_0,x_1]$ and $|\lambda|\leq \varepsilon_1$. We have
$H_2(B)-H_2(A)=O(\mu)$ which gives
$$\displaystyle -\int_{\widehat{AB}}
 F_2(x)G_2'(x)dt=\mu\varphi(a,\lambda,\mu),$$
$$ \begin{array}{ll}
\displaystyle\int_{\widehat{AB}}x^{2j+1}G_2'(x)dt
&=\displaystyle\oint_{\gamma_\lambda(a)}x^{2j+1}G_2'(x)dt +\mu\psi_j
(a,\lambda,\mu),\\
&=-\displaystyle\oint_{\gamma(a)}x^{2j+1}dy
+\lambda\varphi_j(a,\lambda,\mu)+\mu\psi_j
(a,\lambda,\mu),\end{array}
$$ where $\varphi$, $\varphi_j$ and $\psi_j$ are smooth functions
for $a\in[x_0,x_1],$ $|\lambda|\leq \varepsilon_1$ and $|\mu|\leq
\varepsilon_2$. Thus, we have
\begin{equation}\label{27}
H_2(B)-H_2(A)=\mu[M(h)+\lambda\varphi+\displaystyle\sum_{j=0}^qb_j(\lambda\varphi_j+\mu\psi_j)]\equiv\mu
d(h,\lambda,\mu),\end{equation} where
$$\begin{array}{c}
M(h)=\displaystyle\sum_{j=0}^qb_jI_j(h), \ \
I_j(h)=\displaystyle\oint_{H_2=h}x^{2j+1}dy,\ h=G_2(a)\in[h_0,h_1].
\end{array} $$
 Let $\deg G_2(x)=2l$. Then by Lemma 2.3, we see
that
$$
\frac{I_{j+1}}{I_j}=\frac{\alpha_{j+1}}{\alpha_j}h^{\frac{1}{l}}(1+o(1)),\
j=0,1,\cdots,q-1$$ for $h\gg 1$. Hence, we first fix $b_0\neq 0$ and
then vary $b_1$, $b_2$, $\cdots$, $b_q$ in turn satisfying
$$|b_q|\ll |b_{q-1}|\ll \cdots \ll |b_0|, \  b_jb_{j+1}<0, \
j=0,\cdots,q-1$$ such that $M(h)$ has $q$ zeros, denoted by $\bar
h_i$, $i=1,\cdots,q$, on the interval $[h^*,+\infty)$, each having
an odd multiplicity.

Now we fix $b_j$ as taken before, and then choose $N$, $h_1$, $x_1$
and $\bar \varepsilon_0$ as follows:
$$\begin{array}{l}
N=\displaystyle\max_{0\leq j\leq q}\{|b_j|\},\\
h_1=\displaystyle\max_{1\leq i\leq q}\{\bar h_i\}+1=G_2(x_1),\\
\bar
\varepsilon_0=\max\{\varepsilon_1(x_1),\varepsilon_2(x_1,N)\}.\end{array}$$

Then, by the above discussion, the function $d(h,\lambda,\mu)$ in
\eqref{27} is well defined for all $h\in [h_0,h_1]$, $|\lambda|\leq
\bar \varepsilon_0$, $|\mu|\leq \bar \varepsilon_0$, and the
function $M(h)$ has $q$ different zeros with odd multiplicity on the
open interval $(h^*,h_1)$. Let
$$U=\{(x,y)| \ h^*\leq H(x,y)\leq h_1\}.$$

Since all of the $q$ zeros of $M$ have an odd multiplicity, there
exists an $\varepsilon_0\in (0,\bar \varepsilon_0)$  such that for
all $|\lambda|\leq  \varepsilon_0$, $|\mu|\leq  \varepsilon_0$ the
function $d(h,\lambda,\mu)$ has $q$ zeros in $h\in (h^*,h_1)$, each
having an odd multiplicity, and that the corresponding limit cycles
of \eqref{26} are all located in $U$. This ends the proof.

The main result of this section is the following theorem.

{\bf Theorem 2.1.} \  {\it If \eqref{21} has the property $Z(n,m,k)$,
then there exist two polynomial systems  of the form
\begin{equation} \label{28}
\dot x=y-\varepsilon F_{2n+2}, \ \ \dot y=-g_{2m+1}(x),
\end{equation}
and
\begin{equation} \label{281}
\dot x=y-\varepsilon F_{2n+3}, \ \ \dot y=-g_{2m+1}(x),
\end{equation}
which have properties $Z(2n+1,2m+1,2k+n)$ and $Z(2n+2,2m+1,2k+n+1)$,
respectively, where
$$\begin{array}{c}
2n+1\leq \deg F_{2n+2}\leq 2n+2,\ \  \deg F_{2n+3}=2n+3, \\ 3\leq
\deg g_{2m+1}\leq 2m+1, \ \  g_{2m+1}(-x)=-
g_{2m+1}(x).\end{array}$$}

{\bf Proof.}  First, by Lemma 2.2,  there are a constant
$\varepsilon_0>0$ and a compact set $U_0\subset \mathbb{R}^2$ such that for
all $0<\varepsilon\leq \varepsilon_0$ the system \eqref{21} has at
least $2k$ limit cycles in $U_0$, each having an odd multiplicity.
Let
$$F_2(x)= F(x^2+x_0),\ \ G_2(x)=\int_0^x2xg(x^2+x_0)dx.$$
Then $\deg F_2(x)\leq 2n+2$, $\deg G_2(x)=2l$, $2\leq l\leq m+1$.
Consider \eqref{26} with $q=n$ or $n+1$. For each fixed $\lambda$
with $0<|\lambda|\leq \varepsilon_0$, \eqref{26} has $2k$ limit
cycles in $U_0$ for all $0<|\mu|\ll|\lambda|$, each having an odd
multiplicity. On the other hand, by Lemma 2.4, there exist
$\varepsilon^*\in (0,\varepsilon_0)$, constants $b_j\neq 0$,
$j=0,1,\cdots,q$ and a compact set $U$ with $U\bigcap U_0=\emptyset$
such that for all $|\lambda|\leq \varepsilon^*$, $0<|\mu|\leq
\varepsilon^*$, \eqref{26} has $q$ limit cycles in $U$, each having
an odd multiplicity. Therefore, for all $(\lambda,\mu)$ satisfying
\begin{equation}\label{29}
0<|\mu|\ll|\lambda|\leq \varepsilon^*
\end{equation}
\eqref{26} has $2k+q$ limit cycles in the set $U\bigcup U_0$. Then
set $g_{2m+1}(x)=G_2'(x)$, and
$$F_{2n+2}(x)=\lambda F_2(x)+\mu
\sum_{j=0}^nb_jx^{2j+1},\ F_{2n+3}(x)=\lambda F_2(x)+\mu
\sum_{j=0}^{n+1}b_jx^{2j+1}.
$$
It follows that for all $0<\varepsilon\leq 1$ and some
$(\lambda,\mu)$ satisfying \eqref{29}, the system \eqref{28} has
the property $Z(2n+1,2m+1,2k+n)$, and \eqref{281} has property
$Z(2n+2,2m+1,2k+n+1)$. This finishes the proof.

An obvious corollary of the theorem above is the following.

{\bf Corollary 2.1.} \  {\it If \eqref{21} has the property $Z(n,m,k)$,
then there exist two  polynomial systems of the form
$$
\dot x=y-\varepsilon F_{2n+2}, \ \ \dot y=-g_{2m+2}(x)
$$
and $$ \dot x=y-\varepsilon F_{2n+3}, \ \ \dot y=-g_{2m+2}(x)$$
which has properties $Z(2n+1,2m+2,2k+n)$ and $Z(2n+2,2m+2,2k+n+1)$,
respectively.}

The theorem above together with Corollary 2.1 is fundamental. We can
use them repeatedly. See the next two sections.

\section{Estimate of $H(n,m)$ for fixed $m$}

Suppose there is a system of the form \eqref{11} which has property
$Z(n_0, m_0,k_0)$ with $n_0\geq 1$, $m_0\geq 1$. Then define
\begin{equation}\label{2100}\begin{array}{c}
n_{11}=2n_0+1,\ n_{12}=2n_0+2,\\
 m_{11}=2m_0+1,\ m_{12}=2m_0+2, \\
k_{11}=2k_0+n_0, \ k_{12}=2k_0+n_0+1.\end{array}\end{equation} Then
by Theorem 2.1 and Corollary 2.1, there are polynomial Li\'enard
systems of the form \eqref{11} which have 4 properties
$Z(n_{1i},m_{1j},k_{1i})$ for $i,j=1,2$,  respectively, which imply
\begin{equation}\label{21000}H(n_{1i},m_{1j})\geq k_{1i}, \
i,j=1,2.\end{equation}

Hence, using Lemma 2.1 and formulas in \eqref{2100} and
\eqref{21000} we obtain the following.

{\bf Theorem 3.1.} \ {\it We have

$(1)$ \ $H(n,4)\geq H(n,3)\geq 2[\frac{n-1}{4}]+[\frac{n-1}{2}]$,
$n\geq 3;$

$(2)$ \ $H(n,6)\geq H(n,5)\geq 2[\frac{n-1}{3}]+[\frac{n-1}{2}]$,
$n\geq 3.$}

{\bf Proof.} We only prove the first conclusion. The second one can
be shown similarly. For any integer $\tilde n\geq 3$, let
$n=[\frac{\tilde n-1}{2}]$. Then $n\geq 1$, and either $\tilde
n=2n+1$ or $\tilde n=2n+2$. By Lemma 2.1, we have property
$Z(n,1,[\frac{n}{2}])$. Then by \eqref{2100}  we  further obtain
properties $Z(2n+1,j,2[\frac{n}{2}]+n)$ and
$Z(2n+2,j,2[\frac{n}{2}]+n+1)$, $j=3,4$. Using  \eqref{21000} we
conclude
$$H(2n+1,j)\geq 2[\frac{n}{2}]+n,\ H(2n+2,j)\geq 2[\frac{n}{2}]+n+1,
\ j=3,4.$$ It follows that
$$H(\tilde n,j)\geq \displaystyle2[\frac{1}{2}[\frac{\tilde n-1}{2}]]+[\frac{\tilde
n-1}{2}]=2[\frac{\tilde n-1}{4}]+[\frac{\tilde n-1}{2}].$$ This
finishes the proof.

By Theorem 2.1 and Corollary 2.1, we can get more results. For the
purpose, let $$ S_{m_0}=\{2^i(m_0+1)-2+j\,|\,1\leq j\leq 2^i,\,i\geq
0\}
$$ for $ m_0\geq 1$,
It is easy to see that $S_1\bigcup S_2=\{m|\ m\geq 1\}$. Thus by  Lemma 3.2 in \cite{hl},  for any integer $M\geq 1$,
\begin{equation} \label{213}
\bigcup_{m_0=M}^{2M}S_{m_0}=\{m\,|\,m\geq M\}.
\end{equation}

{\bf Theorem 3.2.} \ {\it We have

$(1)$ For any integers $p\geq 1$ and $1\leq j\leq 2^p$ there exist
constants $0<\delta_p<(p+4)2^{p-1}-(p+1)$  and
$0<\beta_p<(p+\frac{13}{3})2^{p-1}-(p+\frac{4}{3})$ such that
\begin{equation}\label{36}H(n,2^{p+1}-2+j)\geq \frac{1}{2}(p+1)n-\delta_p,\end{equation}
\begin{equation}\label{37}H(n,3\cdot2^{p}-2+j)\geq
(\frac{1}{2}p+\frac{2}{3})n-\beta_p.\end{equation} In particular,
\begin{equation}\label{371} H(n,7)\geq\frac{3}{2}n-9 \ \ {\rm for}\ \ n\geq 7. \end{equation}
 $(2)$ For any
integer $m\geq 7$ there exists $\gamma_m>0$ satisfying
$$\lim_{m\rightarrow\infty} \sup \frac{\gamma_m}{(m+1)\ln(m+1)}\leq \frac{1}{4\ln2}$$
such that
\begin{equation}\label{38}H(n,m)\geq
(\frac{\ln(m+2)}{2\ln2}-\frac{1}{3})n-\gamma_m, \ \ n\geq
m.\end{equation} In particular,
$$\displaystyle\lim_{m\rightarrow\infty}\lim_{n\rightarrow\infty} \inf \frac{H(n,m)}{n\ln(m+2)}\geq
\frac{1}{2\ln2}.$$

}

{\bf Proof.} Note that $[\frac{n}{2}]\geq \frac{n}{2}-\frac{1}{2}$
for all $n\geq 1$. We have $$2[\frac{n}{2}]+n\geq
2(\frac{n}{2}-\frac{1}{2})+n=2n-1,\ 2[\frac{n-1}{2}]-1\geq
2(\frac{n-1}{2}-\frac{1}{2})-1=n-3.$$ Hence, the first conclusion of
Theorem 3.1 implies
$$H(n, 2^{2}-2+j)\geq l_1n-\delta_1, \ n\geq 3, j=1,2,$$
where $l_1=1$, $\delta_1=3$. Just following the idea in the proof of
Theorem 3.1, we can obtain further
$$H(2n+2,2^3-1)\geq H(2n+1,2^3-1)\geq 2(l_1n-\delta_1)+n,\ \
n\geq 2^2-1,$$ which implies
$$\begin{array}{cl}
H(n, 2^{3}-1)&\!\!\geq (2l_1+1)[\frac{n-1}{2}]-2\delta_1\\
&\!\!\geq (2l_1+1)(\frac{n-1}{2}-\frac{1}{2})-2\delta_1\\
&\!\!=l_2n-\delta_2, \ n\geq 2^3-1,\end{array}$$ where $l_2=l_1+ \frac{1}{2}=\frac{3}{2}$, $\delta_2=9$. In particular,  \eqref{371} follows.
 Note that $H(n, 2^{3}-2+j)\geq H(n, 2^{3}-1)$ for $1\leq j\leq
2^2.$ Hence, we have
$$H(n, 2^{3}-2+j)\geq l_2n-\delta_2, \ n\geq 2^3-1, \ 1\leq j\leq 2^2.$$
 In the same way, we have for $p \geq 3$
\begin{equation}\label{39}H(n, 2^{p+1}-2+j)\geq l_pn-\delta_p, \ n\geq 2^{p+1}-1, \ 1\leq j\leq 2^p,\end{equation}
where $l_p$ and $\delta_p$ satisfy
$$l_p=l_{p-1}+\frac{1}{2},\ \
\delta_p=2\delta_{p-1}+2l_{p-1}+1.$$Then using these relations and
the initial data $l_1=1$ and $\delta_1=3$, we can easily find
$$l_p=\frac{p+1}{2}, \ \  \delta_p=3\cdot2^{p-1}+\sum_{j=2}^{p}(j+1)2^{p-j}, \ \ p\geq 1.$$
 Then noting
\begin{equation}\label{390}\delta_p\leq3\cdot2^{p-1}+(p+1)\sum_{j=2}^{p}2^{p-j}=(p+4)2^{p-1}-(p+1),\end{equation}
\eqref{36}  follows from \eqref{39}.

Further, using $[\frac{n-1}{3}]\geq\frac{n-1}{3}-\frac{2}{3}$ it
follows from the second conclusion of Theorem 3.1
$$H(n,6)\geq H(n,5)\geq
2(\frac{n-1}{3}-\frac{2}{3})+\frac{n-1}{2}-\frac{1}{2}=r_1n-\beta_1,$$
where $n\geq 5$, $r_1=\frac{7}{6}$, $\beta_1=3$.

Similar to the above,  we can obtain for $p \geq 1$
\begin{equation}\label{310}H(n, 3\cdot2^{p}-2+j)\geq r_pn-\beta_p, \ n\geq 3\cdot2^{p}-1, \ 1\leq j\leq 2^p,\end{equation}
where $r_p$ and $\beta_p$ satisfy
$$r_p=r_{p-1}+\frac{1}{2},\ \
\beta_p=2\beta_{p-1}+2r_{p-1}+1,$$ which together with
$r_1=\frac{7}{6}$, $\beta_1=3$ give
$$r_p=\frac{p}{2}+\frac{2}{3}, \ \ \beta_p=3\cdot2^{p-1}+\sum_{j=2}^{p}(j+\frac{4}{3})2^{p-j}, \ \ p\geq 1.$$
Then noting
$$\beta_p\leq3\cdot2^{p-1}+(p+\frac{4}{3})\sum_{j=2}^{p}2^{p-j}=(p+\frac{13}{3})2^{p-1}-(p+\frac{4}{3}),$$
\eqref{37} follows.

By \eqref{213}, any positive integer $m$ is either in $S_1$ or in
$S_2$. If it is in $S_1$, then by the definition of $S_{m_0}$ there
exist $p\geq 1$ and $1\leq j\leq 2^p$ such that $m=2^{p+1}-2+j$,
which implies $2^{p+1}-1\leq m\leq 2^{p+1}-2+2^p=3\cdot 2^p-2$, or
$p+1\leq \frac{1}{\ln2}\ln(m+1)$, $p\geq
\frac{1}{\ln2}\ln\frac{m+2}{3}$. Hence, by \eqref{39} and
\eqref{390} we have
$$H(n,m)\geq \bar l_mn-\bar \delta_m,\ \ n\geq m, $$
where $$\bar l_m=\frac{1}{2}(1+\frac{1}{\ln2}\ln\frac{m+2}{3}), \ \
\displaystyle\lim_{m\rightarrow\infty}\sup\frac{\bar
\delta_m}{(m+1)\ln(m+1)}\leq \frac{1}{4\ln2}.$$

 If  $m\in S_2$, then similarly we have $m=3\cdot 2^p-2+j$ for some
 $ p\geq 0$ and $1\leq j\leq 2^p$. Thus, $\frac{1}{\ln2}\ln\frac{m+1}{3}\geq p\geq
 \frac{\ln(m+2)}{\ln2}-2$. By \eqref{310} we have
$$H(n,m)\geq \bar r_mn-\bar \beta_m,\ \ n\geq m, $$
where $$\bar r_m=\frac{\ln(m+2)}{2\ln2}-\frac{1}{3}, \ \
\displaystyle\lim_{m\rightarrow\infty}\sup\frac{\bar
\beta_m}{(m+1)\ln(m+1)}\leq \frac{1}{6\ln2}.$$ Therefore, for $n\geq
m$

$$\displaystyle H(n,m)\geq \min\{\bar l_m,\bar r_m\}n-\max\{\bar \delta_m,\bar \beta_m\}$$ which
yields \eqref{38} since $\bar r_m<\bar l_m$. The proof is completed.

In \eqref{36} and \eqref{37}, taking $j=1$ we have in particular

$$\displaystyle H(n,m)\geq \frac{\ln(m+1)}{2\ln2}n-\bar \delta_m,\ \ n\geq m$$ for $m=2^{p+1}-1$
and
$$\displaystyle H(n,m)\geq (\frac{\ln(m+1)}{2\ln2}+\frac{1}{6\ln2}\ln
\frac{16}{9})n-\bar \beta_m, \ \ n\geq m$$ for $m=3\cdot2^{p}-1.$

\section{Estimate of  $H(m,m)$}

Following \eqref{2100}, define further
\begin{equation}\label{2110}\begin{array}{c} n_{21}=2n_{11}+1,\
n_{22}=2n_{11}+2,n_{23}=2n_{12}+1,\
n_{24}=2n_{12}+2,\\
m_{21}=2m_{11}+1,\ m_{22}=2m_{11}+2, \ m_{23}=2m_{12}+1,\
m_{24}=2m_{12}+2,\\
 k_{21}=2k_{11}+n_{11},\ k_{22}=2k_{11}+n_{11}+1,\ k_{23}=2k_{12}+n_{12},
 \ k_{24}=2k_{12}+n_{12}+1.\end{array}\end{equation}
Then by Theorem 2.1 and Corollary 2.1 again, there are polynomial
Li\'enard systems of the form \eqref{11} which have $4^2$ properties
$Z(n_{2i},m_{2j},k_{2i})$, $i,j=1,2,3,4=2^2$, respectively, which
give $$H(n_{2i},m_{2j})\geq k_{2i}, \ i,j=1,2,3,4.$$ Also, we have
obviously from \eqref{2100} and  \eqref{2110}
$$\begin{array}{c}k_{12}=k_{11}+1, \ \ k_{24}= k_{23}+1,\ k_{22}= k_{21}+1,\\
k_{21}=2^2k_0+2^2n_0+1, \ \ k_{23}=2^2k_0+2^2n_0+4. \end{array}
$$

In general, we introduce three  series $n_{pj}$, $m_{pj}$ and
$k_{pj}$ for $p\geq 1$ and $j=1,\cdots, 2^p$. We do it by induction
as follows:
\begin{equation}\label{210}
\begin{array}{c}
n_{i+1,2l-1}=2n_{il}+1, \ \ n_{i+1,2l}=2n_{il}+2,\\
m_{i+1,2l-1}=2m_{il}+1, \ \ m_{i+1,2l}=2m_{il}+2,\\
k_{i+1,2l-1}=2k_{il}+n_{il},
\ \ k_{i+1,2l}=2k_{il}+n_{il}+1,\\
l=1,\cdots,2^i, \ \ i\geq 1.
\end{array}
\end{equation}

Then we can prove the following result.

{\bf Theorem 4.1.} \ {\it If there exists a polynomial system of
 the form \eqref{21} which has the property $Z(n_0,m_0,k_0)$, then for
 all
$p\geq 1$, $i,j=1,\cdots,2^p$
 there are  polynomial systems of
 the form \eqref{21} which have properties respectively $Z(n_{pi},m_{pj},k_{pi})$ for all $p\geq 1$, $i,j=1,\cdots,2^p$, and
 therefore, $H(n_{pi},m_{pj})\geq k_{pi}$,
 where
\begin{equation}\label{211}
\begin{array}{c}
n_{pi}=2^p(n_0+1)-2+i,\ \ m_{pj}=2^p(m_0+1)-2+j,
\end{array}
\end{equation} and
\begin{equation}\label{0211}
\begin{array}{c}
 k_{p1}=2^p(k_0-1)+p2^{p-1}(n_0+1)+1,\\
 k_{p,2^p}=2^p(k_0-1)+p2^{p-1}(n_0+2)+1,\\
k_{p1}< k_{p2}<\cdots <k_{p,2^p}.
\end{array}
\end{equation}

 }

{\bf Proof.} The formulas for $n_{pi}$ and $m_{pj}$ in \eqref{211} follow from  Lemma 3.1 in \cite{hl}. From \eqref{2110} and \eqref{210} it is easy to prove by induction  that
$k_{pj}\geq k_{pi}$ for $1\leq i< j\leq 2^p$. Hence, to finish the proof, it suffices to prove the two equalities in \eqref{0211}.

For the first equality, it is true for $p=1$. Suppose it is true for
$p=i$. Then by \eqref{210} and the formula for $n_{i1}$ we have
$$k_{i+1,1}=2k_{i,1}+n_{i1}=2k_{i,1}+2^i(n_0+1)-1$$
which, together with the inductive assumption, yields that
$$\begin{array}{ll}
k_{i+1,1}\!\!\!&=2(2^i(k_0-1)+i2^{i-1}(n_0+1)+1)+2^i(n_0+1)-1\\
&=2^{i+1}(k_0-1)+(i+1)2^{i}(n_0+1)+1.
\end{array}
$$ This shows that the equality is also true for $p=i+1$. Then the
first equality in \eqref{0211} follows.

The second one can be obtained in the same way since by \eqref{210}
and \eqref{211}
$$k_{i+1,2^{i+1}}=2k_{i,2^i}+n_{i,2^i}+1=2k_{i,2^i}+2^i(n_0+2)-1.$$
This completes the proof.

In the following we take $n_0=m_0$ and suppose that  there exists a
polynomial system of
 the form \eqref{21} which has the property $Z(m_0,m_0,k_0)$. Then by
Theorem 4.1, we have
\begin{equation}\label{41}H(m_{pi},m_{pi})\geq 2^p(k_0-1)+p(m_0+1)2^{p-1}+1,\end{equation}
 where
$m_{pi}=2^p(m_0+1)-2+i,$ $ p\geq 1,$ $1\leq i \leq 2^p$.

Let $m>m_0$ and $m\in S_{m_0}$. Then $m=m_{pi}$ for some $ p\geq 1,$
$1\leq i \leq 2^p$. If $i=1$, then $m=2^p(m_0+1)-1$, or
$$ 2^p=\frac{m+1}{m_0+1}, \ \
p=\frac{1}{\ln2}\ln\frac{m+1}{m_0+1}.$$

Then by \eqref{41} we obtain
\begin{equation}\label{42}H(m,m)\geq\frac{(m+1)\ln(m+1)}{2\ln2}+N(m_0,k_0)(m+1)+1,\end{equation}
for $m=2^p(m_0+1)-1$, $ p\geq 1$, where $$
N(m_0,k_0)=\frac{k_0-1}{m_0+1}-\frac{\ln(m_0+1)}{2\ln2}.$$

If $2\leq i\leq 2^p$, then
$$ 2^p\geq \frac{m+2}{m_0+2}, \ \
p\geq\frac{1}{\ln2}\ln\frac{m+2}{m_0+2}.$$ By \eqref{41} again we
obtain
\begin{equation}\label{43}H(m,m)\geq N_1(m_0)\frac{(m+2)\ln(m+2)}{2\ln2}+N_2(m_0,k_0)(m+2)+1\equiv H_{m_0}(m),\end{equation}
where $$m\in S_{m_0},\ \ N_1(m_0)=\frac{m_0+1}{m_0+2},\ \
N_2(m_0,k_0)=\frac{k_0-1}{m_0+2}-\frac{N_1(m_0)\ln(m_0+2)}{2\ln2}.$$

Similarly, by Theorem 4.1, we have $$H(m_{p,2^p},m_{p,2^p})\geq
2^p(k_0-1)+p(m_0+2)2^{p-1}+1,$$
 where
$m_{p,2^p}=2^p(m_0+2)-2$, which gives
\begin{equation}\label{40} H(m,m)\geq\frac{(m+2)\ln(m+2)}{2\ln2}+N(m_0+1,k_0)(m+2)+1,\end{equation}
for $m=2^p(m_0+2)-2,\ \ p\geq 1$.

By \eqref{42} and  \eqref{40} we can obtain

{\bf Theorem 4.2.} \ {\it For any $k\geq 2$ there exist constants
$B_k$ and $\bar B_k$ satisfying
$$\lim_{k\rightarrow\infty}\frac{B_k}{\ln k}=\lim_{k\rightarrow\infty}\frac{\bar B_k}{\ln
k}=-\frac{1}{2\ln2}$$ such that
$$H(m,m)\geq\frac{(m+1)\ln(m+1)}{2\ln2}+B_k(m+1)+1$$
for $m=2^pk-1$, $ p\geq 1$, and
$$H(m,m)\geq\frac{(m+2)\ln(m+2)}{2\ln2}+\bar B_k(m+2)+1$$
for $m=2^p(k+1)-2$, $ p\geq 1$.

In particular, for $m=2^{p+1}-1$, $p\geq 1$, we have
\begin{equation}\label{44}H(m,m)\geq \frac{(m+1)\ln(m+1)}{2\ln2}+1.\end{equation}
}

{\bf Proof.} We need only to prove \eqref{44}. In fact, by Lemma
2.1, we have the property $Z(3,3,5)$. Thus, taking $k_0=5$, $m_0=3$ in
\eqref{42} we obtain  \eqref{44} directly. The proof is ended.

 By \eqref{43}, we have further

{\bf Theorem 4.3.} \ {\it $(1)$ For all $m\geq 3$ we have

\begin{equation}\label{45}H(m,m)\geq \frac{(m+2)\ln(m+2)}{3\ln2}-\frac{m+2}{3}(1+\frac{\ln3}{\ln2})+1.\end{equation}

$(2)$ \begin{equation}\label{46}
\lim_{m\rightarrow\infty}\inf\frac{H(m,m)}{(m+2)\ln(m+2)}\geq
\frac{1}{2\ln2}.\end{equation} That is to say, $H(m,m)$ grows at
least as rapidly as $\frac{1}{2\ln2}(m+2)\ln(m+2)$ as $m$ goes to
infinity.}

{\bf Proof.} Using properties $Z(1,1,0)$ and $Z(2,2,1)$, by
\eqref{43} we have $H(m,m)\geq H_1(m)$ for $m\in S_1$, $H(m,m)\geq
H_2(m)$ for $m\in S_2$, where
$$\begin{array}{cl}
 H_1(m)&=\displaystyle\frac{(m+2)\ln(m+2)}{3\ln2}-\frac{m+2}{3}(1+\frac{\ln3}{\ln2})+1\\
 &\leq\displaystyle\frac{3(m+2)\ln(m+2)}{8\ln2}-\frac{m+2}{3}(1+\frac{\ln3}{\ln2})+1\\
 &\leq\displaystyle\frac{3(m+2)\ln(m+2)}{8\ln2}-\frac{3(m+2)}{4}+1=H_2(m).
\end{array}$$
It follows that for $m\in S_1\bigcup S_2$ we have $H(m,m)\geq
H_1(m)$, which gives \eqref{45}.

Now let $M$ be an arbitrary integer. Let
$$\bar N_1(M)=\min\{N_1(m_0)|\ M\leq m_0\leq 2M\},$$
$$\bar N_2(M)=\min\{N_2(m_0,0)|\ M\leq m_0\leq 2M\}.$$
Then by \eqref{43}, we have $\bar N_1(M)=\frac{M+1}{M+2}$ and
$$H(m,m)\geq \bar N_1(M)\frac{(m+2)\ln(m+2)}{2\ln2}+\bar
N_2(M)(m+2)+1$$ for all $\displaystyle m\in
\bigcup_{m_0=M}^{2M}S_{m_0}$. By \eqref{213}, the above inequality
holds for all $m\geq M$. Therefore,

$$\lim_{m\rightarrow\infty}\inf\frac{H(m,m)}{(m+2)\ln(m+2)}\geq
\frac{\bar N_1(M)}{2\ln2}.$$  Since $M$ is
arbitrary \eqref{46} follows. This ends the proof.

\section{Estimate of $H(m\pm r,m)$ }

Let $r$ be a positive integer. We can give an estimate of $H(m\pm
r,m)$ using the result obtained in the previous section or the  method
used above. First, noting
 $$H(m+r,m)\geq H(m,m), \ \ H(m-r,m)\geq H(m-r,m-r),$$ we obtain by \eqref{44} and
 \eqref{46}
$$H(m+r,m)\geq \frac{(m+1)\ln(m+1)}{2\ln2}+1$$
for $m=2^{p+1}-1$, $p\geq 1$;
$$H(m-r,m)\geq \frac{(m+1-r)\ln(m+1-r)}{2\ln2}+1$$
for $m=2^{p+1}-1+r$, $p\geq 1$; and
\begin{equation}\label{50}\lim_{m\rightarrow\infty}\inf\frac{H(m\pm r,m)}{m\ln m}\geq
\frac{1}{2\ln2}.\end{equation}

However, for the case $H(m-r,m)$ using the above method, we can get
more and better estimate. We describe the process briefly here.

First, suppose that there exists a polynomial system of
 the form \eqref{21} which has the property $Z(m_0-r,m_0,k_0)$.
 Then, as we do in Theorem 2.1, we can construct a system of the  form \eqref{26} which has the property $Z(2q,2m_0+1,2k_0+q)$. It follows that $H(2q,2m_0+1)\geq 2k_0+q$.
In order to obtain $H(2m_0+1-r,2m_0+1)\geq 2k_0+q$, we need to have
$2q\leq 2m_0+1-r$ or $q\leq m_0-[\frac{r}{2}]$. We take $q=
m_0-[\frac{r}{2}]$ and introduce
$$ m_{11}=2m_0+1,\ m_{12}=2m_0+2, \
k_{11}=2k_0+m_0-[\frac{r}{2}]=k_{12}.$$ Then we have properties
$Z(m_{1j}-r,m_{1j},k_{1j})$, $j=1,2$. In general, define
$$\begin{array}{c}
m_{i+1,2l-1}=2m_{il}+1, \ \ m_{i+1,2l}=2m_{il}+2,\\
k_{i+1,2l-1}=k_{i+1,2l}=2k_{il}+m_{il}-[\frac{r}{2}],\\
l=1,\cdots,2^i, \ \ i\geq 1.\end{array} $$

Then, as before, we have $H(m_{ij}-r,m_{ij})\geq k_{ij}$  where
\begin{equation}\label{51}
\begin{array}{c}
m_{ij}=2^i(m_0+1)-2+j, \ i\geq 1, \  j=1,\cdots,2^i, \\
 k_{i1}=2^i(k_0-1-[\frac{r}{2}])+i2^{i-1}(m_0+1)+1+[\frac{r}{2}],\\
 k_{i,2^i}=2^i(k_0-2-[\frac{r}{2}])+i2^{i-1}(m_0+2)+2+[\frac{r}{2}],\\
k_{i1}\leq k_{i2}\leq \cdots \leq k_{i,2^i}.
\end{array}
\end{equation}

Therefore, using \eqref{51}, similar to \eqref{42}, \eqref{40} and
 \eqref{43} we can obtain
\begin{equation}\label{52}H(m-r,m)\geq\frac{(m+1)\ln(m+1)}{2\ln2}+N(m_0,k_0-[\frac{r}{2}])(m+1)+1+[\frac{r}{2}]\end{equation}
for $m=2^i(m_0+1)-1$, $ i\geq 1$, and
\begin{equation}\label{53}H(m-r,m)\geq\frac{(m+2)\ln(m+2)}{2\ln2}+N(m_0+1,k_0-1-[\frac{r}{2}])(m+2)+2+[\frac{r}{2}]\end{equation}
for $m=2^i(m_0+2)-2$, $ i\geq 1$, and
\begin{equation}\label{54}H(m-r,m)\geq N_1(m_0)\frac{(m+2)\ln(m+2)}{2\ln2}+N_2(m_0,k_0-[\frac{r}{2}])(m+2)+1+[\frac{r}{2}]\equiv \bar H_{m_0}(m)\end{equation}
for $m\in S_{m_0}$.

 Using \eqref{52}, \eqref{53} and \eqref{54},
just similar to Theorems 4.2 and 4.3 we can obtain the following theorems.

{\bf Theorem 5.1.} \ {\it Let $r$ be a positive integer. For any $k\geq 2$ there exist constants $B_{k,r}$ and $\bar B_{k,r}$ satisfying
$$ B_{k,r}\geq-\Big[\frac{1}{k}\Big(1+[\frac{r}{2}]\Big)+\frac{\ln k}{2\ln2}\Big],$$
$$\bar B_{k,r}\geq-\Big[\frac{1}{k+1}\Big(2+[\frac{r}{2}]\Big)+\frac{\ln(k+1)}{2\ln2}\Big],$$
such that
$$H(m-r,m)\geq\frac{(m+1)\ln(m+1)}{2\ln2}+B_{k,r}(m+1)+1+[\frac{r}{2}]$$
for $m=2^pk-1$, $ p\geq 1$, and
$$H(m-r,m)\geq\frac{(m+2)\ln(m+2)}{2\ln2}+\bar B_{k,r}(m+2)+2+[\frac{r}{2}]$$
for $m=2^p(k+1)-2$, $ p\geq 1$. }

{\bf Theorem 5.2.} \ {\it $(1)$ For $m=2^{p+1}-1$, $p\geq 1$, we
have $$H(m-1,m)\geq \frac{(m+1)\ln(m+1)}{2\ln2}+1.$$

$(2)$  For all $m\geq 3$ we have $$H(m-1,m)\geq
\frac{(m+2)\ln(m+2)}{3\ln2}-\frac{m+2}{3}(1+\frac{\ln3}{\ln2})+1.$$
}

We mention that in Theorem 5.2 we have taken  $r=1$. In this case,
$H(m-1,m)$ denotes the maximal number of limit cycles of polynomial
Li\'enard systems of degree $m$.


\end{document}